\documentclass[a4paper,draft]{article}
\usepackage{a4}
\usepackage{amsfonts}
\usepackage{amsmath}
\usepackage{amssymb}

\newcommand{\heute}{9 September 2003}

\numberwithin{equation}{section}
\numberwithin{figure}{section}

\newtheorem{theorem}[equation]{Theorem}
\newtheorem{lemma}[equation]{Lemma}
\newtheorem{corollary}[equation]{Corollary}
\newtheorem{proposition}[equation]{Proposition}

\newenvironment{proof}[1][]{\begin{trivlist} \item[\hskip\labelsep
\indent \emph{Proof#1.}]}{\foorp \end{trivlist}}
\newcommand{\foorp}{{\unskip\nobreak\hfil\penalty50
 \hskip1em\vadjust{}\nobreak\hfil \vrule height5pt width5pt depth0pt
 \parfillskip=0pt \finalhyphendemerits=0 \par}}

\newenvironment{notation}{\begin{trivlist} \item[\hskip\labelsep
\indent \emph{Notation.}]}{\end{trivlist}}

\newenvironment{oneofftheorem}[1]{\begin{trivlist} \item[\hskip\labelsep
\textbf{#1}]}{\end{trivlist}}

\newcommand{\f}[1][p]{\mathbb{F}_{#1}}

\newcommand{\Yal}[1]{Yal\-\c{c}{\i}n}
\newcommand{\Poi}[1]{Poincar\'{e}}
\newcommand{\Gro}[1]{Gr\"ob\-ner}
\newcommand{\diag}{\operatorname{diag}}
\newcommand{\Ess}{\operatorname{Ess}}
\newcommand{\Tr}{\operatorname{Tr}}

\newcommand{\SU}{\mathit{SU}}

\newenvironment{textmatrix}{\left(\begin{smallmatrix}}{\end{smallmatrix}\right)}
\newcommand{\Sq}{\operatorname{Sq}}
\newcommand{\Einf}{E_{\infty}}
\newcommand{\N}{\mathcal{N}}

\begin{document}

\author{David J. Green}
\title{The essential ideal in group cohomology does not square to zero}
\date{\heute}

\maketitle

\begin{abstract}
\noindent
Let $G$ be the Sylow $2$-subgroup of the unitary group $\SU_3(4)$.
We find two essential classes in the mod~$2$ cohomology ring of~$G$
whose product is nonzero. In fact, the product is the ``last survivor''
of Benson--Carlson duality. Recent work of Pakianathan and \Yal. then
implies a result about connected graphs with an action of~$G$.
Also, there exist essential classes which
cannot be written as sums of transfers from proper subgroups.

This phenomenon was first observed on the computer.
The argument given here uses the elegant calculation by J.~Clark,
with minor corrections.
\end{abstract}

\section{Introduction}
Let $G$ be a finite group and $k$ a field of characteristic~$p$. A class
$x$~in the cohomology ring $H^*(G,k)$ is called essential if its restriction
to every proper subgroup is zero. Standard properties of the transfer mean that
$G$ has to be a $p$-group if there are to be nonzero essential classes.

The essential classes in $H^*(G,k)$ form an ideal, called the essential ideal
$\Ess(G)$. If the group~$G$ is elementary abelian, then there is a
non-nilpotent essential class: the product of the Bocksteins
of all nonzero elements of $H^1$.
But in all remaining cases that could be calculated,
the pattern that emerged seemed to be that every product of essential classes
was zero. This led H.~Mui and T.~Marx to formulate the following

\begin{oneofftheorem}{Essential Conjecture} (\cite{Mui:Essay,Marx:Ess2})
Let $G$ be a  finite $p$-group which is not elementary abelian.
Then $\Ess(G)^2 = 0$.
\end{oneofftheorem}

The evidence for the Essential Conjecture consists of a wealth of worked
examples, together with one result which is valid for all finite groups:
Minh showed in~\cite{Minh:Ess-p} that if $x$~is essential then $x^p = 0$
(assuming $G$~is not elementary abelian).

The Essential Conjecture recently acquired new significance through the
work of Pakianathan and \Yal.~\cite{PakYal:NilpotentIdeals} on
$G$-CW-complexes, as we shall see in Corollary~\ref{coroll:PakYal}.
\medskip

\noindent
In this paper we disprove the Essential
Conjecture. That is, we prove the following:

\begin{theorem}
\label{theorem:main}
Let $G$ be a Sylow $2$-subgroup of the special unitary group $\SU_3(4)$.
Then $\Ess(G)^2 \not = \{0\}$. To be more precise, there are essential
classes in $H^4(G,\f[2])$ and $H^{10}(G,\f[2])$ whose product is nonzero.
This nonzero element of $H^{14}(G,\f[2])$ is the last survivor in the sense
of Benson and Carlson~\cite{BensonCarlson:Poincare}.
Moreover there are essential classes in degrees six and eight whose product
is the last survivor too.
\end{theorem}

\noindent
Theorem~\ref{theorem:main} is proved in Propositions \ref{prop:Ess10:4}~and
\ref{prop:Ess8:6}.
Two corollaries can be immediately deduced.
The first one was predicted by Minh in~\cite{Minh:Ess-p}, who observed that
in cases studied to date, the essential ideal had to square to zero, because
each essential class was a sum of transfers from proper subgroups. This raises
the question of whether essential classes are always sums of transfers.

\begin{corollary}
Let $G$ be a Sylow $2$-subgroup of the special unitary group $\SU_3(4)$.
Then $\Ess(G)$ is not contained in $\Tr(G)$, the ideal in
$H^*(G,\f[2])$ generated by all transfers from proper subgroups.
\end{corollary}

\begin{proof}
By Frobenius reciprocity, essential classes annihilate transfer classes.
So if two essential classes have nonzero product, then neither can lie in
$\Tr(G)$.
\end{proof}

\noindent
The second corollary is particularly intriguing. It concerns the theory
of transformation groups and was predicted by Pakianathan and \Yal.
in~\cite{PakYal:NilpotentIdeals}. Recall from tom~Dieck's
book \cite[III (3.14)]{tomDieck:book} that if $G$ is elementary abelian and $X$
a finite-dimensional $G$-CW-complex, then $X$~has a fixed point if and only if
the map $H_G^*({*}) \rightarrow H_G^*(X)$ is injective. The key assumption
here is that $G$ is elementary abelian. Now consider the following corollary:

\begin{corollary}
\label{coroll:PakYal}
Let $G$ be a Sylow $2$-subgroup of the special unitary group $\SU_3(4)$
and let $X$ be a connected $1$-dimensional $G$-CW-complex. Then $X$ has
a fixed point if and only if the map $H_G^*({*}) \rightarrow H_G^*(X)$
is injective.
\end{corollary}

\begin{proof}
This is Theorem~6.2 of~\cite{PakYal:NilpotentIdeals}.
\end{proof}

\noindent
This raises the interesting question of whether there are other groups for
which $\Ess(G)^r$ is nonzero for values of $r$ greater than two. For the groups
concerned this would lead by the results of~\cite{PakYal:NilpotentIdeals}
to similar fixed point theorems in higher dimensions. An obvious family
of groups to study would be the Sylow $2$-subgroups of $\SU_3(2^n)$.
\medskip

\noindent
The fact that $\Ess(G)^2$ is nonzero for the Sylow $2$-subgroups of $\SU_3(4)$
was first observed on the computer~\cite{habil}. The argument presented here
uses J.~Clark's elegant spectral sequence calculation~\cite{Clark:U34}, with
minor corrections.

\section{Clark's approach to the cohomology calculation}
\label{section:ClarksApproach}
Let $G$ be a Sylow $2$-subgroup of the unitary group $\SU_3(4)$.
See Aschbacher's book \cite[\S22]{Aschbacher:book} for a discussion of unitary
groups, and Evens' book \cite{Evens:book} for a general reference on group
cohomology.  As in Exercise~1 on p.~100 of~\cite{Aschbacher:book},
Clark identifies $G$ with the group of upper triangular matrices
$\begin{textmatrix} 1 & a & b \\ 0 & 1 & a^4 \\ 0 & 0 & 1 \end{textmatrix}$
over $\f[16]$ satisfying $b + b^4 = a^5$.
The centre of~$G$ is elementary abelian of rank two, the quotient
by the centre is elementary abelian of rank four, and every noncentral
element has exponent four. So the centre is the sole maximal elementary
abelian subgroup and the cohomology ring over any field of characteristic~$2$
is Cohen--Macaulay
\label{page:Duflot}
by Duflot's theorem (Theorem 10.3.1 in Evens' book~\cite{Evens:book}).

The cohomology ring of $G$ over the field $\f[16]$ was calculated
by Clark~\cite{Clark:U34} using the Lyndon--Hochschild--Serre spectral
sequence of the central extension $1 \rightarrow Z(G) \rightarrow G
\rightarrow G/Z(G) \rightarrow 1$. Clark makes elegant use of the Frobenius
automorphism and of the action of the normalizer of~$G$.

The normalizer of~$G$ in $\SU_3(4)$ is $WG$,
where $W$ is the cyclic group
of order~$15$ generated by the diagonal matrix $T :=
\diag(\zeta^{-1}, \zeta^{-3}, \zeta^4)$ for $\zeta \in \f[16]$ is a primitive
$15$th root of unity. Choose $\zeta$ to be a zero of $X^4 + X^3 + 1$,
meaning that
\begin{equation}
\label{eqn:Nzeta}
\zeta^8 + \zeta^4 + \zeta^2 + \zeta = 1 \, .
\end{equation}
The spectral sequence is equivariant for the conjugation action of~$W$,
and the action of $T$~on
$H^1(G,\f[16]) = H^1(G/Z(G),\f[16])$ is diagonalisable with
eigenvalues $\zeta, \zeta^2, \zeta^4, \zeta^8$.
Similarly the action of $T$~on $H^1(Z(G),\f[16])$
is diagonalisable with eigenvalues $\zeta^5, \zeta^{10}$.
The Frobenius automorphism $F \colon \f[16] \rightarrow \f[16]$,
$F(\alpha) = \alpha^2$ also acts on the spectral sequence,
and the spectral sequence
for $H^*(G,\f[16])$ is just the spectral sequence for $H^*(G,\f[2])$ tensored
with $\f[16]$. Moreover, the actions of $T$~and $F$ commute with each other.

Pick a $T$-eigenvector $a_1 \in H^1(G, \f[16])$ with eigenvalue~$\zeta$,
and set $a_2 = F(a_1)$, $a_4 = F(a_2)$ and $a_8 = F(a_4)$. Then
$a_1,a_2,a_4,a_8$ is a basis for $H^1(G,\f[16])$ which is permuted cyclically
by~$F$. Here and elsewhere in this paper, a subscript~$i$ means that the
element concerned is an eigenvector for~$T$ with eigenvalue $\zeta^i$.

Pick $u_5 \in H^1(Z(G),\f[16])$ to be an eigenvector of~$T$ with
eigenvalue~$\zeta^5$. Then $u_5 \in E_2^{0,1}$ and $d_2(u_5)$
must be a scalar multiple of~$a_1 a_4$, as this is a basis for the eigenspace
in $E_2^{2,0}$ with eigenvalue~$\zeta^5$. As $G$~is not elementary abelian,
$d_2(u_5)$ cannot be zero, and so we may choose $u_5$ to satisfy
$d_2(u_5) = a_1 a_4$. Note that $u_5 \in H^1(Z(G),\f[4])$, since
$d_2(F^2(u_5)) = d_2(u_5)$. Setting $u_{10} = F(u_5)$, we have
that $u_5,u_{10}$ is a basis%
\footnote{%
We write $u_5,u_{10}$ for Clark's $a_5,a_{10}$. Apart from this we
follow his notation.}
for $H^1(Z(G),\f[16])$ which is transposed by~$F$. Moreover,
\begin{equation}
\label{eqn:d2}
d_2(u_5) = a_1 a_4 \quad \text{and} \quad d_2(u_{10}) = a_2 a_8 \, .
\end{equation}
An element $x \in E^{s,t}_r$ is rational (meaning: comes from
the spectral sequence for~$\f[2]$) if and only if $F(x) = x$.
So in particular
\begin{equation}
\label{eqn:H1F2}
H^1(G,\f[2]) =
\{ \alpha a_1 + \alpha^2 a_2 + \alpha^4 a_4 + \alpha^8 a_8 \mid
\alpha \in \f[16] \} \, .
\end{equation}

\begin{notation}
Let us denote by $\N$ the $\f[2]$-linear map
$\N(x) = x + F(x) + F^2(x) + F^3(x)$. Then $\N \colon H^n(G,\f[16]) \rightarrow
H^n(G,\f[2])$ and $\N \colon E_r^{pq}(\f[16]) \rightarrow E_r^{pq}(\f[2])$,
the latter being a map of spectral sequences which is $H^*(G,\f[2])$-linear
but does not in general respect products.
\end{notation}


\section{Essential classes}
\label{section:essential}
We shall now assume that the reader is familiar with Clark's spectral sequence
calculation of $H^*(G,\f[16])$, at least as far as the $\Einf$-page. 
This calculation is performed in Appendix~\ref{section:ClarkCorrected} for
the sake of completeness,
and because we need an explicit basis for the $\Einf$-page.
The appendix also points out some errors in Clark's relations.

Observe that the ideal of essential classes in $H^*(G,\f[16])$ is obtained
from the essential ideal in $H^*(G,\f[2])$ by extending scalars. So for
the Essential Conjecture it is immaterial whether we work over a prime field
or over an extension field.

A sufficient (and in fact necessary) condition for a class in $H^*(G,\f[16])$
to be essential is that it is divisible by each nonzero
element of $H^1(G,\f[2])$. So a class is essential if it is divisible by
$\N(\alpha a_1)$ for each $\alpha \in \f[16]^*$
(c.f.\@ Eqn.~\eqref{eqn:H1F2}).

\begin{lemma}
\label{lemma:H4Ess}
Every inflation class in $H^4(G,k)$ is essential.
\end{lemma}

\begin{proof}
A basis for $E_{\infty}^{4,0}$ consists of
$a_1^4,a_2^4,a_3^4,a_4^4,a_1^3 a_2,a_2^3 a_4, a_4^3 a_8, a_8^3 a_1$.
Now,
\[
a_1^2 a_2 \N(\alpha a_1) = \alpha a_1^3 a_2
\qquad
a_1^3 \N(\alpha a_1) = \alpha a_1^4  + \alpha^2 a_1^3 a_2
\]
by the relations in Equations \eqref{eqn:E4Rels}~and \eqref{eqn:E4MoreRels}.
So $a_1^3 a_2$~and $a_1^4$ are both divisible by
$\N(\alpha a_1)$ for each $\alpha \in \f[16]^*$.
Now apply the Frobenius map~$F$.
\end{proof}

\begin{notation}
We shall say that a class in $E_{\infty}^{**}$ is $\Einf$-essential
if it is divisible by each nonzero $F$-stable element of $E^{1,0}_{\infty}$. 
\end{notation}

Let $x \in \Einf^{pq}$ be an $\Einf$-essential class. For certain
values of $p,q$ we shall show that every class $y \in H^*(G,\f[16])$ which
is represented by~$x$ is essential.
For $q=0$ this is obvious.

\begin{lemma}
\label{lemma:E64}
Each class in $E_{\infty}^{6,2}$ and in $\Einf^{6,4}$ is $E_{\infty}$-essential.
Hence each element of $H^8(G,\f[16])$ which is represented by an element
of $\Einf^{6,2}$ is essential, as is each element of $H^{10}(G,\f[16])$
which is represented by an element of $E_{\infty}^{6,4}$.
\end{lemma}

\begin{proof}
The second part follows from the first because
$\Einf^{6,2}$ and $E_{\infty}^{6,4}$ are subspaces of $H^8(G,\f[16])$ and
$H^{10}(G,\f[16])$ respectively
(see Figure~\ref{fig:EinfTable} on p.~\pageref{fig:EinfTable})\@.
For the first part, note that the multiplication maps
$\Einf^{2,2} \otimes \Einf^{4,0} \rightarrow \Einf^{6,2}$ and
$\Einf^{2,4} \otimes \Einf^{4,0} \rightarrow \Einf^{6,4}$ are surjective.
But every element of $\Einf^{4,0}$ is $\Einf$-essential,
by the proof of Lemma~\ref{lemma:H4Ess}.
\end{proof}

\begin{corollary}
\label{coroll:EinfEssential}
Suppose that $x \in \Einf^{4,4}$ is $\Einf$-essential. Then every class
$y \in H^8(G,\f[16])$ which is represented by~$x$ is essential. Similarly,
if $x \in \Einf^{6,4}$ is $\Einf$-essential, then every $y \in H^{10}(G,k)$
represented by~$x$ is essential.
\end{corollary}

\begin{proof}
Let $z \in H^1(G,\f[2]) \setminus \{0\}$.
By assumption there is $w \in \Einf^{3,4}$ such that
$zw = x$. Let $\tilde{w} \in H^7(G,\f[16])$ be a class which is represented
by~$w$. Then $y - z\tilde{w}$ lies in the subspace
$\Einf^{6,2}$~of $H^8(G,\f[16])$.
So by Lemma~\ref{lemma:E64} there is $v \in H^7(G,\f[16])$ with
$y = z(\tilde{w}+v)$. The second part is proved similarly.
\end{proof}

\begin{lemma}
\label{lemma:H8Ess}
Every element of $\Einf^{4,4}$ is $\Einf$-essential.
Therefore every element of $H^8(G,\f[16])$,
whose restriction to the centre of~$G$ is trivial, is essential.
\end{lemma}

\begin{proof}
The $F$-orbits of $a_4^2 \delta_3$ and $a_4^2 \delta_7$ constitute
a basis for $\Einf^{4,4}$.  One has
\[
a_4 \delta_3 \N(\alpha a_1)
= \alpha^4 a_4^2 \delta_3 + \alpha^8 a_4^2 \delta_7
\quad \text{and} \quad
a_4 \delta_7 \N(\alpha a_1)
= \alpha^4 a_4^2 \delta_7
\]
from the definitions in Eqn.~\eqref{eqn:EinfGenDefs} and
the relations in the $E_4$-page. Now use the corollary.
\end{proof}

\noindent
Note that $\N$~satisfies the product
rule $\N(A) \cdot \N(B) = \N(A \cdot \N(B))$.

\begin{lemma}
\label{lemma:H10Ess}
Suppose that $y \in H^{10}(G,\f[16])$ is represented by a class
$x \in \Einf^{4,6}$ of the form $x = \N(\omega a_4 \tau_6)$
with $\omega \in \f[4]^*$.  Then $y$ is an essential class.
\end{lemma}

\begin{proof}
The class~$x$ is $\Einf$-essential since
$\N(\omega \alpha^{11} \tau_6) \cdot \N(\alpha a_1)
= x + \N(\omega \alpha^4 a_8 \tau_6)$ and
$(\omega \tau_{10} + \omega^2 \tau_5) \N(\alpha a_1)
= \N(\omega \alpha^4 a_8 \tau_6)$.
Now use Corollary~\ref{coroll:EinfEssential}.
\end{proof}

\noindent
Observe that $(X^8+X^4+X^2+X)(X^8+X^4+X^2+X+1)=X^{16}+X$ in $\f[2][X]$,
and that therefore $X^8 + X^4+X^2+X$ has eight distinct zeros in~$\f[16]$.
Observe further that these zeros form an $\f[2]$-vector subspace.

\begin{lemma}
\label{lemma:H6Ess}
Suppose that $y \in H^6(G,\f[16])$ is represented by a class $x \in \Einf^{4,2}$
of the form
$x = \N(\lambda a_4^2 \beta_7)$ with $\lambda \in \f[16]$ satisfying
$\N(\lambda) = 0$. Then $y$~is essential.
\end{lemma}

\begin{proof}
Since $\Einf^{4,2}$ is a subspace of $H^6$, it suffices to show that
$x$~is $\Einf$-essential.
Observe that $\prod(X^2+X+\alpha) = X^{16} + X$ in $\f[2][X]$, where the
product is taken over the zeros $\alpha$~of $X^8 + X^4 + X^2 + X$. So
we may pick a $\mu \in \f[16]$ satisfying $\mu^2 + \mu = \lambda$.
One then has
$\N(\alpha a_1) \cdot \N(\lambda \alpha^{11} a_4 \beta_7 + \mu a_8 \beta_7)
= \N(\lambda a_4^2 \beta_7)$ for all $\alpha \in \f[16]^*$.
\end{proof}

\begin{proposition}
\label{prop:Ess10:4}
Pick $\omega \in \f[4]^*$.
Let $\eta \in H^{10}(G,\f[2])$ be the class $\N(\eta')$
for some $\eta' \in H^{10}(G,\f[16])$ which is represented by
$\omega a_4 \tau_6 \in E_{\infty}^{4,6}$.
Let $\theta \in H^4(G,\f[2])$ be the class
$\N(\omega^{-1} \zeta a_4^3 a_8)$.
Then $\eta$~and $\theta$ are both essential classes, and $\eta \theta = \xi$.
Hence $\Ess(G)^2$ is nonzero.
\end{proposition}

\noindent
See \S\ref{subsection:LastSurvivor} for the definition of the class~$\xi$.
It is the last survivor in the sense
of Benson and Carlson~\cite{BensonCarlson:Poincare} for the polynomial
generators representing $u_5^8,u_{10}^8$.

\begin{proof}
Lemma~\ref{lemma:H4Ess} says that $\theta$~is essential, and $\eta$~is
essential by Lemma~\ref{lemma:H10Ess}.
Note that
$\omega a_4 \tau_6 \N(\omega^{-1} \zeta a_4^3 a_8)
= \zeta a_4^5 a_8 \beta_7 u_{10}^4 = \zeta \xi$,
and so $\eta \theta = \N(\zeta \xi)$.
As $\xi$~is $F$-stable and $\N(\zeta) = 1$ by Eqn.~\eqref{eqn:Nzeta},
it follows that $\N(\zeta \xi) = \xi$.
\end{proof}

\begin{proposition}
\label{prop:Ess8:6}
Pick $\lambda \in \f[16]^*$ to be one of the seven zeros of $X^7+X^3+X+1$.
Let $\phi \in H^6(G,\f[2])$ be $\N(y)$ for $y = \lambda a_4^2 \beta_7$
in $\Einf^{4,2} \subseteq H^6(G,\f[16])$.
Let $\psi \in H^8(G, \f[2])$ be $\N(z)$ for some $z \in H^8(G,\f[16])$
which is represented by
$\lambda^{-1} \zeta a_4^2 \delta_7 \in \Einf^{4,4}$.
Then $\phi$~and $\psi$ are essential, and $\phi \psi = \xi$.
\end{proposition}

\begin{proof}
Lemma~\ref{lemma:H6Ess} says that $\phi$~is essential, and $\psi$~is essential
by Lemma~\ref{lemma:H8Ess}.
Now $\lambda a_4^2 \beta_7 \N(\lambda^{-1} \zeta a_4^3 a_8 u_{10}^4)
= \zeta a_4^4 \tau_6 = \zeta \xi$. Therefore
$\N(y) \N(z) = \N(\zeta) \xi = \xi$.
\end{proof}

\paragraph{Report on the computer calculation}
The fact that $\Ess(G)^2$ is nonzero for this group~$G$ was first
observed using the author's computer program~\cite{habil}.
As a module over the polynomial algebra $k[u_5^8,u_{10}^8]$, the ideal of
essential classes is free with 75 generators. The degrees of these generators
are given by the generating function
$8t^4 + 6t^5 + 3t^6 + 8t^7 + 16t^8 + 7t^9 + 6t^{10} + 8t^{11} + 8t^{12}
+ 4t^{13} + t^{14}$.
The calculations above show that in dimensions
4,5,6,8,10 and 14 there are (at least) as many essential classes as the
computer states there to be.
The computer states further that Propositions \ref{prop:Ess10:4}~and
\ref{prop:Ess8:6} account for all nontrivial products of essential classes.
In particular, $\Ess(G)^3$ is zero.

\appendix

\section{Clark's spectral sequence calculation}
\label{section:ClarkCorrected}
For Section~\ref{section:essential} we need an explicit basis for the
$\Einf$-page of the spectral sequence. We shall derive such a basis here.
In \S\ref{subsection:ClarkErrors} we point out
some errors in Clark's original calculation~\cite{Clark:U34}.

From Section~\ref{section:ClarksApproach} we know that
the $E_2$-page of the spectral sequence
is the polynomial algebra $\f[16][a_1,a_2,a_4,a_8,u_5,u_{10}]$, with the
first four generators lying in $E_2^{1,0}$ and the last two in $E_2^{0,1}$.
From the formula for the differential~$d_2$ in Eqn.~\eqref{eqn:d2} one sees
that the $E_3$-page is $\f[16][a_1,a_2,a_4,a_8]/(a_1 a_4, a_2 a_8)
\otimes \f[16][u_5^2,u_{10}^2]$: for
$a_1 a_4, a_2 a_8$ is a regular sequence in~$E_2^{**}$.

\subsection{The $E_4$-page of the spectral sequence}
\label{subsect:E4}
\begin{lemma}
\label{lemma:E4}
The $E_4$-page of the spectral sequence is generated by
$a_1,a_2,a_4,a_8 \in E_4^{1,0}$, $\beta_7,\beta_{14},\beta_{13},\beta_{11}
\in E_4^{2,2}$ and $u_5^4,u_{10}^4 \in E_4^{0,4}$. These generators are permuted
by~$F$ in the obvious way. The classes $\beta_i$ are represented in $Z_5^{2,2}$
by
\begin{equation}
\label{eqn:BetaDef}
\beta_7 = a_4 a_8 u_5^2 \quad \beta_{14} = a_8 a_1 u_{10}^2
\quad
\beta_{13} = a_1 a_2 u_5^2 \quad \beta_{11} = a_2 a_4 u_{10}^2
\, .
\end{equation}
The relations ideal is generated by
\begin{equation}
\label{eqn:E4Rels}
\begin{array}{c}
a_1 a_4 \quad a_2 a_8
\qquad a_8 a_1^2 + a_2 a_4^2 \quad a_4 a_8^2 + a_1 a_2^2 \\
a_1 \beta_7 \quad a_2 \beta_7 \quad a_2 \beta_{14} \quad a_4 \beta_{14}
\quad a_4 \beta_{13} \quad a_8 \beta_{13}
\quad a_8 \beta_{11} \quad a_1 \beta_{11} \\
a_1 \beta_{14} + a_2 \beta_{13} + a_4 \beta_{11} + a_8 \beta_7
\qquad \beta_i^2 \quad \beta_i \beta_j
\end{array}
\; .
\end{equation}
The following relations also hold:
\begin{equation}
\label{eqn:E4MoreRels}
\begin{array}{c}
a_1^2 a_2^2 \quad a_2^2 a_4^2 \quad a_4^2 a_8^2 \quad a_8^2 a_1^2
\qquad
a_1 a_2^3 \quad a_2 a_4^3 \quad a_4 a_8^3 \quad a_8 a_1^3
\\
a_8^2 \beta_7 + a_8 a_1 \beta_{14}
\quad
a_1^2 \beta_{14} + a_1 a_2 \beta_{13}
\quad
a_2^2 \beta_{13} + a_2 a_4 \beta_{11}
\quad
a_4^2 \beta_{11} + a_4 a_8 \beta_7
\\
a_4 a_8^2 \beta_7 \qquad a_8 a_1^2 \beta_{14} \qquad a_1 a_2^2 \beta_{13}
\qquad a_2 a_4^2 \beta_{11}
\end{array} \, .
\end{equation}
The $E_4$-page is $E_4^{0*} \otimes (E_4^{*0} \oplus E_4^{*2})$.
A basis for $E_4^{*0}$ is
\begin{equation}
\label{eqn:E4star0Basis}
\begin{array}{c}
a_1^{r+1}, \, a_2^{r+1}, \, a_4^{r+1}, \, a_8^{r+1}; \;
a_1^{r+1} a_2, \, a_2^{r+1} a_4, \, a_4^{r+1} a_8, \, a_8^{r+1} a_1
\quad \text{for $r \geq 0$;} \\
1; \qquad a_1 a_2^2, \, a_2 a_4^2
\end{array} \, .
\end{equation}
This basis is permuted by~$F$. One orbit has length one, one has length two,
and the remainder have length four.
A basis for $E_4^{*2}$ is
\begin{equation}
\label{eqn:E4star2Basis}
\left. \begin{array}{l}
a_4^r \beta_7, \, a_8^r \beta_{14}, \, a_1^r \beta_{13}, \, a_2^r \beta_{11} \\
a_4^{r+1} a_8 \beta_7, \, a_8^{r+1} a_1 \beta_{14},
\, a_1^{r+1} a_2 \beta_{13}, \, a_2^{r+1} a_4 \beta_{11}
\end{array} \right\}
\, \text{for $r \geq 0$;}
\quad a_8 \beta_7, \, a_2 \beta_{13}, \, a_4 \beta_{11} \, .
\end{equation}
The orbit sum of~$a_8 \beta_7$ is zero. The rest of the basis is permuted freely
by~$F$.
\end{lemma}

\begin{proof}
A basis for $E_3^{*0}$ is: $1$; $a_1^{r+1}$~and its $F$-orbit (length four);
$a_1^{r+1} a_2^{s+1}$~and its $F$-orbit (length four). Here $r,s\geq 0$.
Note that $x^2 = F \circ \Sq^r(x)$ for $x \in H^r(G,\f[16])$.  Hence
$d_3(u_5^2) = d_3 \Sq^1(u_{10}) = \Sq^1 d_2 u_{10} = a_1^2 a_8 + a_2 a_4^2$,
and so $d_3$~acts as follows:
\begin{alignat*}{3}
u_5^2 & \mapsto a_8 a_1^2 + a_2 a_4^2
& \qquad
a_1^{r+1} u_5^2 & \mapsto a_8 a_1^{r+3}
& \qquad
a_1^{r+1} a_2^{s+1} u_5^2 & \mapsto 0
\\
u_{10}^2 & \mapsto a_1 a_2^2 + a_4 a_8^2
& \qquad
a_1^{r+1} u_{10}^2 & \mapsto a_1^{r+2} a_2^2
& \qquad
a_1^{r+1} a_2^{s+1} u_{10}^2 & \mapsto a_1^{r+2} a_2^{s+3}
\end{alignat*}
Hence $E_4^{*0}$ has the advertised basis, and a basis for $Z_3^{*2}$
consists of the $a_1^r a_2^s \beta_{13}$ for $r,s\geq 0$, together with
their length four $F$-orbits.
Since $a_i a_j a_k = 0$ in~$E_3^{3,0}$ if $i,j,k$  are pairwise distinct,
one sees that $\beta_i \beta_j = 0$ in~$E_3^{4,4}$ if $i \not = j$; and
$a_i \beta_j = 0$ in $E_3^{3,2}$ if $a_i$~does not feature in the definition
of~$\beta_j$.
Since $\beta_7^2 = d_3(a_4 u_5^4 u_{10}^2)$, one has
$\beta_i^2 = 0$ in $E_4^{4,4}$.
The action of $d_3$ on $E_3^{*4}$ is given by
\[
\begin{array}{c}
u_5^2 u_{10}^2 \mapsto a_1 \beta_{14} + a_2 \beta_{13} + a_4 \beta_{11}
+ a_8 \beta_7 \\
a_1^{r+1} u_5^2 u_{10}^2 \mapsto a_1^{r+2} \beta_{14} + a_1^{r+1} a_2 \beta_{13}
\qquad
a_1^{r+1} a_2^{s+1} u_5^2 u_{10}^2 \mapsto a_1^{r+1} a_2^{s+2} \beta_{13}
\end{array}
\, .
\]
So $E_4^{*2}$ has the advertised basis, $Z_3^{*4} = u_5^4 Z_3^{*0}
+ u_{10}^4 Z_3^{*0}$, and the relations ideal is as claimed.
\end{proof}

\subsection{The $\Einf$-page of the spectral sequence}
\label{subsect:E6}
\begin{proposition}
\label{prop:E6}
The $\Einf$-page of the spectral sequence is generated
by
$u_5^8, u_{10}^8 \in \Einf^{0,8}$,
$a_1,a_2,a_4,a_8 \in \Einf^{1,0}$,
$\beta_7,\beta_{14},\beta_{13},\beta_{11} \in \Einf^{2,2}$,
$\delta_3,\delta_6,\delta_{12},\delta_9$ and
$\delta_7,\delta_{14},\delta_{13},\delta_{11}$ in $\Einf^{2,4}$,
$\tau_5,\tau_{10}$ and $\tau_3,\tau_6,\tau_{12},\tau_9$ in $\Einf^{3,6}$,
and $\chi_5,\chi_{10} \in \Einf^{3,8}$.
These generators are permuted by~$F$ in the obvious way.
Their representatives in $Z_5^{**}$ are as follows:
\begin{equation}
\label{eqn:EinfGenDefs}
\begin{array}{c}
\begin{alignedat}{4}
\delta_3 & = a_4^2 u_{10}^4 & \quad \delta_6 & = a_8^2 u_5^4 & \quad
\delta_{12} & = a_1^2 u_{10}^4 & \quad \delta_9 & = a_2^2 u_5^4 \\
\delta_7 & = a_4 a_8 u_{10}^4 & \quad \delta_{14} & = a_8 a_1 u_5^4 & \quad
\delta_{13} & = a_1 a_2 u_{10}^4 & \quad \delta_{11} & = a_2 a_4 u_5^4 \\
\tau_3  & = a_2 \beta_{11} u_5^4 & \quad \tau_6 & = a_4 \beta_7 u_{10}^4
& \quad
\tau_{12} & = a_8 \beta_{14} u_5^4 & \quad \tau_9 & = a_1 \beta_{13} u_{10}^4
\end{alignedat} \\
\begin{alignedat}{2}
\tau_5 & = (a_8 \beta_7 + a_4 \beta_{11}) u_5^4 & \qquad
\tau_{10} & = (a_2 \beta_{13} + a_4 \beta_{11}) u_{10}^4 \\
\chi_5 & = a_1 a_2^2 u_5^4 u_{10}^4
& \qquad \chi_{10} & = a_2 a_4^2 u_5^4 u_{10}^4
\end{alignedat}
\end{array}
\, .
\end{equation}
The $\Einf$-page is $\Einf^{0*} \otimes (\Einf^{*0} \oplus \Einf^{*2} \oplus \Einf^{*4}
\oplus \Einf^{*6} \oplus \f[16]\{\chi_5,\chi_{10}\})$.
A basis for $\Einf^{*0}$ is
\begin{equation}
\label{eqn:E6star0Basis}
\begin{array}{c}
a_1^{r+1}, \, a_2^{r+1}, \, a_4^{r+1}, \, a_8^{r+1}; \;
a_1^{r+1} a_2, \, a_2^{r+1} a_4, \, a_4^{r+1} a_8, \, a_8^{r+1} a_1
\quad \text{for $r \leq 3$;} \\
1; \quad a_1 a_2^2, \, a_2 a_4^2; \quad a_1^5, a_2^5
\end{array} \, .
\end{equation}
This basis is permuted by~$F$. One orbit has length one, two have length two,
and the remainder have length four.
A basis for $\Einf^{*2}$ is
\begin{equation}
\label{eqn:E6star2Basis}
\begin{array}{rl}
a_4^r \beta_7, \, a_8^r \beta_{14}, \, a_1^r \beta_{13}, \, a_2^r \beta_{11}
& \text{for $r \leq 4$} \\
a_4^{r+1} a_8 \beta_7, \, a_8^{r+1} a_1 \beta_{14},
\, a_1^{r+1} a_2 \beta_{13}, \, a_2^{r+1} a_4 \beta_{11}
& \text{for $r \leq 2$}
\end{array}
\quad
a_8 \beta_7, \, a_2 \beta_{13}, \, a_4 \beta_{11}
\, .
\end{equation}
The orbit sum of~$a_8 \beta_7$ is zero. The rest of the basis is permuted freely
by~$F$.
A basis for $\Einf^{*4}$ is
\begin{equation}
\label{eqn:E6star4Basis}
\begin{array}{c}
a_4^r \delta_7, \, a_8^r \delta_{14}, \, a_1^r \delta_{13}, \, a_2^r \delta_{11}
\; \text{for $r \leq 4$}
\qquad
a_4^r \delta_3, \, a_8^r \delta_6, \, a_1^r \delta_{12}, \, a_2^r \delta_9
\; \text{for $r \leq 2$} \\
a_2 \delta_3, \, a_4 \delta_6 \qquad a_8 \delta_7, \, a_1 \delta_{14}
\qquad a_4^3 \delta_3, \, a_8^3 \delta_6, \, a_2^3 \delta_9
\end{array}
\end{equation}
Apart from the last three terms, this basis is permuted by~$F$.
Two orbits have length two, the rest length four.
The orbit sum of $a_4^3 \delta_3$ is zero.
A basis for $\Einf^{*6}$ is
\begin{equation}
\label{eqn:E6star6Basis}
\begin{array}{rl}
a_2^r \tau_3, \, a_4^r \tau_6, \, a_8^r \tau_{12}, \, a_1^r \tau_9
& \text{for $r \leq 3$} \\
a_2^r a_4 \tau_3, \, a_4^r a_8 \tau_6, \, a_8^r a_1 \tau_{12}, \,
a_1^r a_2 \tau_9 & \text{for $r \leq 3$}
\end{array} \quad \tau_5, \, \tau_{10}; \quad \xi \, ,
\end{equation}
where the class $\xi = a_2^4 a_4 \tau_3$ is $F$-stable, and the remaining
basis elements are permuted in the obvious way.

Summing up, the dimension of $\Einf^{pq}/(u_5^8,u_{10}^8)$ may be read off from
Figure~\ref{fig:EinfTable} (an empty entry means dimension zero).
\begin{figure}
\centering
$\begin{array}{c|c|c|c|c|c|c|c|c|c|c}
q \\
\cline{2-10}
8 & & & & 2 & & & & & & \\
\cline{2-10}
7 & & & & & & & & & & \\
\cline{2-10}
6 & & & & 6 & 8 & 8 & 8 & 4 & 1 & \\
\cline{2-10}
5 & & & & & & & & & & \\
\cline{2-10}
4 & & & 8 & 12 & 8 & 7 & 4 & & & \\
\cline{2-10}
3 & & & & & & & & & & \\
\cline{2-10}
2 & & & 4 & 7 & 8 & 8 & 8 & & & \\
\cline{2-10}
1 & & & & & & & & & & \\
\cline{2-10}
0 & 1  & 4 & 8 & 10 & 8 & 6 & & & & \\
\cline{2-11}
\multicolumn{1}{c}{}
& \multicolumn{1}{c}{0}
& \multicolumn{1}{c}{1}
& \multicolumn{1}{c}{2}
& \multicolumn{1}{c}{3}
& \multicolumn{1}{c}{4}
& \multicolumn{1}{c}{5}
& \multicolumn{1}{c}{6}
& \multicolumn{1}{c}{7}
& \multicolumn{1}{c}{8}
& \multicolumn{1}{c}{p}
\end{array}$
\caption{The dimension of the $\Einf$-page modulo $(u_5^8, u_{10}^8)$}
\label{fig:EinfTable}
\end{figure}

The relations ideal for the $\Einf$-page can be recovered from the
relations in the $E_4$-page
(Equations \eqref{eqn:E4Rels}~and \eqref{eqn:E4MoreRels}), together with
the relations
$a_1^5 + a_4^5$, $a_2^5 + a_8^5$ and
$a_4^3 \delta_3 + a_8^3 \delta_6 + a_1^3 \delta_{12} + a_2^3 \delta_9$:
these three expressions are the images under~$d_5$ of~$u_5^4$, $u_{10}^4$~and
$u_5^4 u_{10}^4$ respectively.
\end{proposition}

\begin{proof}
The $\beta_i$ are permanent cycles for degree reasons, and the classes
$u_5^4,u_{10}^4$ are transgressive.
So the $E_5$-page coincides with the $E_4$-page.
We shall show that the $E_6$-page takes the form claimed for the
$\Einf$-page.
Then as $E_6^{r+6,0} = E_6^{r+7,2} = 0$ for all $r \geq 0$, one deduces that
there are no further differentials and that $E_6^{**} = \Einf^{**}$.

Arguing as for~$d_3$, one has
$d_5(u_5^4) = a_1^5 + a_4^5$.
So $d_5 \colon E_5^{*4} \rightarrow E_5^{*0}$ acts thus:
\begin{alignat*}{3}
u_5^4 & \mapsto a_4^5 + a_1^5
& \qquad a_1^{r+1} a_2 u_5^4 & \mapsto a_1^{r+6} a_2
& \qquad a_1 a_2^2 u_5^4 & \mapsto 0 \\
u_{10}^4 & \mapsto a_8^5 + a_2^5
& \qquad a_1^{r+1} a_2 u_{10}^4 & \mapsto 0
& \qquad a_1 a_2^2 u_{10}^4 & \mapsto 0 \\
a_1^{r+1} u_5^4 & \mapsto a_1^{r+6}
& \qquad a_1 u_{10}^4 & \mapsto a_8^5 a_1
& \qquad a_1^{r+2} u_{10}^4 & \mapsto 0
\end{alignat*}
Hence $E_6^{*0}$ has the advertised basis, and a basis for $Z_5^{*4}$
consists of the $F$-orbits of $a_2 \delta_3, a_8 \delta_7$ (both length two)
and the $F$-orbits of $a_4^r \delta_3, a_4^r \delta_7$ (both length four) for
$r \geq 0$. The map $d_5 \colon E_5^{*8} \rightarrow E_5^{*4}$ is given by
\[
\begin{array}{c}
u_5^4 u_{10}^4 \mapsto a_4^3 \delta_3 + a_8^3 \delta_6 + a_1^3 \delta_{12}
+ a_2^3 \delta_9 \\
\begin{alignedat}{2}
a_4 u_5^4 u_{10}^4 & \mapsto a_4^4 \delta_3 + a_2^4 \delta_{11}
& \qquad a_4^{r+2} u_5^4 u_{10}^4 & \mapsto a_4^{r+5} \delta_3 \\
a_1 a_2^2 u_5^4 u_{10}^4 & \mapsto 0
& \qquad a_4^{r+1} a_8 u_5^4 u_{10}^4 & \mapsto a_4^{r+5} \delta_7
\end{alignedat}
\end{array} \, .
\]
This yields the desired basis for $E_6^{*4}$, and shows that $\chi_5,\chi_{10}$
are a basis for a complement of $u_5^8 E_6^{*0} + u_{10}^8 E_6^{*0}$
in $E_6^{*8}$.

The map $d_5 \colon E_5^{*6} \rightarrow E_5^{*2}$ is given by
\begin{alignat*}{2}
a_4^r \beta_7 u_5^4 & \mapsto a_4^{r+5} \beta_7
& \qquad \beta_7 u_{10}^4 & \mapsto a_8^4 a_1 \beta_{14}
\qquad a_4^{r+1} \beta_7 u_{10}^4 \mapsto 0 \\
a_4^{r+1} a_8 \beta_7 u_5^4 & \mapsto a_4^{r+6} a_8 \beta_7
& \qquad a_4^{r+1} a_8 \beta_7 u_{10}^4 & \mapsto 0 \\
a_8 \beta_7 u_5^4 & \mapsto a_4^5 a_8 \beta_7
& \qquad a_8 \beta_7 u_{10}^4 & \mapsto a_8^5 a_1 \beta_{14}
\end{alignat*}
(For the last statement, note that $a_8^6 \beta_7 = a_8^5 a_1 \beta_{14}$
by Eqn.~\eqref{eqn:E4MoreRels}.)
Hence $E_6^{*2}$ has the advertised basis, and a basis for $Z_5^{*6}$
consists of the $F$-orbits of $\tau_5$ (length two) and of
$a_2^r \tau_3, a_2^r a_4 \tau_3$ (both length four) for $r \geq 0$.
The map $d_5 \colon E_5^{*,10} \rightarrow E_5^{*6}$ is given by
\begin{alignat*}{2}
\beta_7 u_5^4 u_{10}^4 & \mapsto a_4^4 \tau_6 + a_8^3 a_1 \tau_{12}
& \qquad a_4^{r+1} \beta_7 u_5^4 u_{10}^4 & \mapsto a_4^{r+5} \tau_6 \\
a_4^{r+1} a_8 \beta_7 u_5^4 u_{10}^4 & \mapsto a_4^{r+5} a_8 \tau_6
& \qquad a_8 \beta_7 u_5^4 u_{10}^4 & \mapsto a_4^4 a_8 \tau_6
+ a_8^4 a_1 \tau_{12}
\end{alignat*}
So $E_6^{*6}$ has the advertised basis, and $E_6^{*,10} = u_5^8 E_6^{*2}
\oplus u_{10}^8 E_6^{*2}$.
\end{proof}

\subsection{The last survivor}
\label{subsection:LastSurvivor}
Let $k$~be a field of characteristic two.
Since $H^*(G,k)$ is a Cohen--Macaulay ring
(as noted on p.~\pageref{page:Duflot}),
one sees from Figure~\ref{fig:EinfTable} that its \Poi. series $P(t)$ is
\begin{equation}
\label{eqn:Poincare}
\textstyle\frac{1+4t+8t^2+10t^3+12t^4+13t^5+16t^6+
20t^7+16t^8+13t^9+12t^{10}+10t^{11}
+8t^{12}+4t^{13}+t^{14}}{(1-t^8)^2} \, .
\end{equation}
Then $P(t)$ satisfies the functional equation
$P(1/t) = t^2 P(t)$, as required by Theorem~1.1 of Benson and Carlson's
paper~\cite{BensonCarlson:Poincare}. This theorem states additionally that
$H^*(G,k)/(u_5^8,u_{10}^8)$ satisfies \Poi. duality\footnote{If $k$~does not
contain $\f[4]$, then $u_5^8,u_{10}^8$ are not defined over~$k$. Consider
instead the parameter system $u_5^8+u_{10}^8, \omega u_5^8 + \omega^2 u_{10}^8$
for $\omega \in \f[4] \setminus \f[2]$. This is defined over~$\f[2]$.}
in formal dimension~$14$, the dual pairing being induced by the cup product
$H^r(G,k) \otimes H^{14-r}(G,k) \rightarrow H^{14}(G,k)$.

\begin{notation}
Denote by $\xi \in H^{14}(G,\f[2])$ the unique nonzero $F$-stable element
of $\Einf^{8,6} \subseteq H^{14}(G,\f[16])$.
This class is the last survivor of \cite[\S7]{BensonCarlson:Poincare}.
By Theorem~1.3 of that paper, it is a transfer from $Z(G)$.
\end{notation}

This class~$\xi$ is represented by
$a_4^4 \tau_6 = a_4^5 a_8 \beta_7 u_{10}^4$,
and by any term in its $F$-orbit.

\subsection{Errors in Clark's calculation}
\label{subsection:ClarkErrors}
Unfortunately some of the relations in the $E_{\infty}$-page are incorrect
on pp.~1424--5 of Clark's paper~\cite{Clark:U34}. Specifically
there are problems with the linear dependencies among the terms
$\delta_i a_j$~in $E_{\infty}^{3,4}$ and among the terms
$\tau_i a_j$~in $E_{\infty}^{4,6}$.
The correct relations in $\Einf^{3,4}$ are that the following twenty expressions
are zero:
\begin{equation}
\label{eqn:Einf34dependencies}
\begin{array}{c}
a_1 \delta_7 \quad a_2 \delta_7 \quad a_2 \delta_{14} \quad a_4 \delta_{14}
\quad a_4 \delta_{13} \quad a_8 \delta_{13}
\quad a_8 \delta_{11} \quad a_1 \delta_{11} \\
a_1 \delta_3 \quad a_2 \delta_6 \quad a_4 \delta_{12} \quad a_8 \delta_9 \\
a_8 \delta_{12} + a_2 \delta_3 \quad a_1 \delta_9 + a_4 \delta_6
\qquad
a_1 \delta_{14} + a_4 \delta_{11} \quad a_2 \delta_{13} + a_8 \delta_7
 \\
a_8 \delta_3 + a_4 \delta_7 \quad a_1 \delta_6 + a_8 \delta_{14} \quad
a_2 \delta_{12} + a_1 \delta_{13} \quad a_4 \delta_9 + a_2 \delta_{11}
\end{array} \, .
\end{equation}
The correct relations in $\Einf^{4,6}$ are that the following sixteen
expressions are zero:
\begin{equation}
\label{eqn:Einf46dependencies}
\begin{array}{c}
a_1 \tau_3 \quad a_2 \tau_6 \quad a_4 \tau_{12} \quad a_8 \tau_9
\qquad
a_8 \tau_3 \quad a_1 \tau_6 \quad a_2 \tau_{12} \quad a_4 \tau_9 \\
a_1 \tau_5 \quad a_2 \tau_{10} \quad a_4 \tau_5 \quad a_8 \tau_{10} \\
a_2 \tau_5 + a_4 \tau_3 \quad a_4 \tau_{10} + a_8 \tau_6 \quad
a_8 \tau_5 + a_1 \tau_{12} \quad a_1 \tau_{10} + a_2 \tau_9
\end{array} \, .
\end{equation}
These relations may be verified using the definitions of the $\delta_i, \tau_j$
(Eqn.~\eqref{eqn:EinfGenDefs}) and the relations in the $E_4$-page
(Equations \eqref{eqn:E4Rels}~and \eqref{eqn:E4MoreRels}). There can be no
more such relations, for one sees from Figure~\ref{fig:EinfTable} that
$\dim \Einf^{3,4}=12$ and $\dim \Einf^{4,6}=8$.

By contrast, Clark states in Figure~1 on p.~1424 of~\cite{Clark:U34} that
$\Einf^{3,4}$ has dimension ten. This leads to an incorrect coefficient of~$t^7$
in the Poincar\'e series at the top of p.~1426: he has $18$, whereas the correct
value is $20$ (see Eqn.~\eqref{eqn:Poincare}).

I have not investigated the ungrading problem for the corrected relations;
in this paper we just need the $E_{\infty}$-page.

\paragraph{Acknowledgements}
The author thanks J. F. Carlson, who first drew this group to his attention.
The computer calculations mentioned here were performed on a machine
belonging to Carlson.

\bibliographystyle{abbrv}
\bibliography{../united}

\begin{thebibliography}{10}

\bibitem{Aschbacher:book}
M.~Aschbacher.
\newblock {\em Finite group theory}, volume~10 of {\em Cambridge Studies in
  Advanced Mathematics}.
\newblock Cambridge University Press, Cambridge, 1986.

\bibitem{BensonCarlson:Poincare}
D.~J. Benson and J.~F. Carlson.
\newblock Projective resolutions and {P}oincar\'e duality complexes.
\newblock {\em Trans. Amer. Math. Soc.}, 342(2):447--488, 1994.

\bibitem{Clark:U34}
J.~Clark.
\newblock Mod-$2$ cohomology of the group ${U}_3(4)$.
\newblock {\em Comm. Algebra}, 22(4):1419--1434, 1994.

\bibitem{Evens:book}
L.~Evens.
\newblock {\em The cohomology of groups}.
\newblock Oxford Univ.\@ Press, Oxford, 1991.

\bibitem{habil}
D.~J. Green.
\newblock {\em {G}r\"ob\-ner Bases and the Computation of Group Cohomology}.
\newblock Lecture Notes in Mathematics. Springer-Verlag, accepted.

\bibitem{Marx:Ess2}
T.~Marx.
\newblock The restriction map in cohomology of finite $2$-groups.
\newblock {\em J. Pure Appl. Algebra}, 67(1):33--37, 1990.

\bibitem{Minh:Ess-p}
P.~A. Minh.
\newblock Essential mod-$p$ cohomology classes of $p$-groups: an upper bound
  for nilpotency degrees.
\newblock {\em Bull. London Math. Soc.}, 32(3):285--291, 2000.

\bibitem{Mui:Essay}
H.~Mui.
\newblock The mod $p$ cohomology algebra of the extra-special group {$E(p^3)$}.
\newblock Unpublished essay, 1982.

\bibitem{PakYal:NilpotentIdeals}
J.~Pakianathan and E.~Yal{\c{c}}{\i}n.
\newblock On nilpotent ideals in the cohomology ring of a finite group.
\newblock {\em Topology}, 42(5):1155--1183, 2003.

\bibitem{tomDieck:book}
T.~tom Dieck.
\newblock {\em Transformation groups}, volume~8 of {\em de Gruyter Studies in
  Mathematics}.
\newblock Walter de Gruyter \& Co., Berlin, 1987.

\end{thebibliography}

\end{document}